\documentclass[english,12pt]{smfart}
\usepackage{amssymb}
\usepackage{amsfonts}
\usepackage{amscd}
\usepackage[T1]{fontenc}
\usepackage[all]{xypic}

\swapnumbers

\def\limproj{\mathop{\oalign{lim\cr
\hidewidth$\longleftarrow$\hidewidth\cr}}}

\def\ac{{\rm ac}}

\let\got\mathfrak
\def\gP{{\got P}}

\def\CC{{\mathbb C}}

\def\FF{{\mathbb F}}
\def\GG{{\mathbb G}}

\def\NN{{\mathbf N}}

\def\QQ{{\mathbb Q}}

\def\ZZ{{\mathbb Z}}

\def\cA{{\mathcal A}}

\def\cC{{\mathcal C}}

\def\cF{{\mathcal F}}

\def\cH{{\mathcal H}}

\def\cL{{\mathcal L}}

\def\cP{{\mathcal P}}

\def\cX{{\mathcal X}}

\mathchardef\alphag="7C0B
\mathchardef\betag="7C0C
\mathchardef\gammag="7C0D
\mathchardef\deltag="7C0E
\mathchardef\varepsilong="7C22
\mathchardef\varphig="7C27
\mathchardef\psig="7C20
\mathchardef\zetag="7C10
\mathchardef\epsilong="7C0F
\mathchardef\rhog="7C1A
\mathchardef\taug="7C1C
\mathchardef\upsilong="7C1D
\mathchardef\iotag="7C13
\mathchardef\thetag="7C12
\mathchardef\pig="7C19
\mathchardef\sigmag="7C1B
\mathchardef\etag="7C11
\mathchardef\omegag="7C21
\mathchardef\kappag="7C14
\mathchardef\lambdag="7C15
\mathchardef\mug="7C16
\mathchardef\xig="7C18
\mathchardef\chig="7C1F
\mathchardef\nug="7C17
\mathchardef\varthetag="7C23
\mathchardef\varpig="7C24
\mathchardef\varrhog="7C25
\mathchardef\varsigmag="7C26
\mathchardef\Omegag="7C0A
\mathchardef\Thetag="7C02
\mathchardef\Sigmag="7C06
\mathchardef\Deltag="7C01
\mathchardef\Phig="7C08
\mathchardef\Gammag="7C00
\mathchardef\Psig="7C09
\mathchardef\Lambdag="7C03
\mathchardef\Xig="7C04
\mathchardef\Pig="7C05
\mathchardef\Upsilong="7C07

\newtheorem{theorem}[subsubsection]{Theorem}
\newtheorem{lem}[subsubsection]{Lemma}

\newtheorem{prop}[subsubsection]{Proposition}

\newtheorem{monoconjecture}[subsubsection]{Monodromy conjecture}
\newtheorem{gmonoconjecture}[subsubsection]{Generalised monodromy conjecture}
\newtheorem{holoconjecture}[subsubsection]{Holomorphy conjecture}
\newtheorem{gholoconjecture}[subsubsection]{Generalised holomorphy conjecture}

\theoremstyle{definition}

\newtheorem{def-prop}[subsubsection]{Proposition-Definition}
\newtheorem{def-theorem}[subsubsection]{Theorem-Definition}

\theoremstyle{remark}
\newtheorem{remark}[subsubsection]{Remark}

\theoremstyle{plain}

\numberwithin{equation}{subsection}

\def\boxit#1#2{\setbox1=\hbox{\kern#1{#2}\kern#1}%
\dimen1=\ht1 \advance\dimen1 by #1
\dimen2=\dp1 \advance\dimen2 by #1
\setbox1=\hbox{\vrule height\dimen1 depth\dimen2\box1\vrule}%
\setbox1=\vbox{\hrule\box1\hrule}%
\advance\dimen1 by .4pt \ht1=\dimen1
\advance\dimen2 by .4pt \dp1=\dimen2 \box1\relax}

\renewcommand{\theequation}{\thesubsection.\arabic{equation}}

\let\got\mathfrak

\newcommand{\N}{\ensuremath{\mathbb{N}}}
\newcommand{\A}{\ensuremath{\mathbb{A}}}
\newcommand{\Z}{\ensuremath{\mathbb{Z}}}
\newcommand{\G}{\ensuremath{\mathbb{G}}}
\newcommand{\Q}{\ensuremath{\mathbb{Q}}}

\newcommand{\C}{\ensuremath{\mathbb{C}}}

\DeclareMathOperator*{\Spec}{Spec}

\begin{document}
\title[Zeta functions and Alexander modules]{Zeta functions and Alexander modules}
\author{Johannes Nicaise}
\address{KULeuven, Dept. of Mathematics, Celestijnenlaan 200B, 3001 Leuven, Belgium}
\email{johannes.nicaise@wis.kuleuven.ac.be}

\begin{abstract}
We introduce the \'etale framework to study Igusa zeta functions
in several variables, generalizing the machinery of vanishing
cycles in the univariate case. We define the \'etale Alexander
modules, associated to a morphism of varieties $F:X\rightarrow
\A^{r}_{k}$, a geometric point $\xi$ of $X$, and an object in the
derived category
$D^{b}_{c}(F^{-1}(\GG^{r}_{m,k}),\bar{\mathbb{Q}}_{l})$. These are
sheaves of modules on the scheme $\mathcal{C}(\GG^{r}_{m,k})$ of
continuous characters $\pi_{1}(\GG^{r}_{m,k})^{t}\rightarrow
\bar{\mathbb{Q}}_l^{\times}$. We formulate Loeser's Monodromy and
Holomorphy Conjectures for multivariate $p$-adic zeta functions,
and prove them in the case where $dim(X)=2$, generalizing results
from the univariate case. Furthermore, we prove a comparison
theorem with the transcendent case, we study a formula of Denef's
for the zeta function in terms of a simultaneous embedded
resolution, and we generalize a result concerning the degree of
the zeta function to our setting.
\end{abstract}

\maketitle

\section{Introduction}
Let us recall the classical definition of a (univariate) local
zeta function, as it was given by Igusa
\cite{Igusa1}\cite{Igusa2}\cite{Igusa3}\cite{Igusa:intro}. A
survey of the univariate theory can be found in Denef's Bourbaki
report \cite{DenefBour}.

In its simplest form, a local zeta function is defined by the
$p$-adic integral
$$Z(f,s)=\int_{\Z_p^{n}}|f(x)|^{s}|dx|\,, $$
where $s$ is a complex variable, $f$ is a polynomial over $\Q_p$
in $n$ variables, $|f|$ is its $p$-adic norm, and $|dx|$ denotes
the Haar measure on the compact group $\Z_p^{n}$, normalized
to give $\Z_p^{n}$ measure $1$. A priori,
$Z(f,s)$ is only defined when $\Re(s)>0$. However, Igusa proved,
using resolution of singularities, that it has a meromorphic
continuation to the complex plane. Moreover, it is a rational
function in $p^{-s}$. The numerical data of an embedded resolution
for $f$ yield a set of candidate poles of $Z(f,s)$, but, since
this set depends on the chosen resolution, a lot of these
candidate poles will not be actual poles of $Z(f,s)$.

When $f$ is defined over $\Q$, Igusa's Monodromy Conjecture
\ref{monoconjecture} predicts an interesting connection between
the eigenvalues of the monodromy at some point of
$f^{-1}(0)\subset \C^{n}$, and the poles of $Z(f,s)$, for almost
all primes $p$, where $f$ is considered as a complex, respectively
$p$-adic polynomial. When the zeta function is twisted by
introducing a multiplicative character
$\chi:\Z_p^{\times}\rightarrow \C^{\times}$, Denef's Holomorphy
Conjecture \ref{holoconjecture} states that poles only occur when
the order of $\chi$ divides the order of an eigenvalue of the
monodromy at some point of $f^{-1}(0)\subset \C^{n}$.

In this paper, we consider multivariate zeta functions, associated
to an $r$-tuple of polynomials $F=(f_1,\ldots,f_r)$. They are
defined in Section \ref{def}. Multivariate local zeta functions
were already studied by Loeser in \cite{ENS}. We state his
generalizations of the Monodromy Conjecture \ref{gmonoconjecture}
and the Holomorphy Conjecture \ref{gholoconjecture}. Since in
general, $F$ does not induce a locally trivial fibration near a
singular point of the complex zero locus, we no longer have a
geometric Milnor fiber at our disposal. However, the formalism of
nearby cycles still makes sense in the multivariate setting, as
was shown by Sabbah in \cite{Sabbah}, such that we do have
something like a cohomological Milnor fiber: the Alexander complex
(actually, it already exists as a motivic object, see
\cite{Guibert}). The monodromy action makes these cohomological
objects into $\C[\Z^{r}]$-modules, and the support of the induced
coherent sheaves on $Spec\,\C[\Z^{r}]$ generalizes the eigenvalues
of the monodromy transformation. Sabbah also proved a generalized
A'Campo formula, expressing the monodromy zeta function of the
Alexander complex at a singular point in terms of a simultaneous
embedded resolution for the polynomials $f_i$.

Using ideas from \cite{GabLoe}, we introduce the \'etale
counterpart of this Alexander complex. The Alexander modules are
associated to a morphism of varieties $F:X\rightarrow \A^{r}_{k}$,
a geometric point $\xi$ of $X$, and an object in the derived
category $D^{b}_{c}(F^{-1}(\GG^{r}_{m,k}),\bar{\mathbb{Q}}_{l})$.
They are sheaves of modules on the scheme
$\mathcal{C}(\GG^{r}_{m,k})$ of continuous characters
$\pi_{1}(\GG^{r}_{m,k})^{t}\rightarrow
\bar{\mathbb{Q}}_l^{\times}$.
 We prove a comparison
theorem, stating that our definition is equivalent to the
 definition of the \'etale nearby cycles if $r=1$, and to
the transcendent one if our base field is the complex field $\C$.
In the latter case, the stalk of the $q$-th Alexander module at a
continuous character $\chi:\pi_{1}(\GG^{r}_{m,k})^{t}\rightarrow
\bar{\mathbb{Q}}_l^{\times}$ corresponds to the subspace of the
$q$-th transcendent Alexander module on which the semi-simple part
of the monodromy acts as multiplication by $\chi$.

Let us give an overview of the structure of the paper. In Section
\ref{def}, we define the multivariate local zeta functions, and we
recall a result of Igusa's concerning their rationality, and the
set of candidate poles associated to a simultaneous embedded
resolution (Theorem \ref{rat}). In Section \ref{alex}, we define
the Alexander modules, and we prove that their support is
contained in a finite union of translated cotori of codimension
one in the scheme of continuous characters
$\mathcal{C}(\G^{r}_{m,k})$ (Proposition \ref{support}). The
generalized Monodromy and Holomorphy Conjectures are stated in
Section \ref{conj}, while Section \ref{trans} contains the
comparison theorems mentioned above (Proposition \ref{comp1} and
\ref{monoco}), making use of a local computation in the normal
crossings case. In Section \ref{form}, we study a generalization
of Denef's formula for the local zeta function in terms of the
geometry of the exceptional locus of a simultaneous embedded
resolution (Theorem \ref{formdenef}). This allows us to prove a
criterion to unmask certain candidate poles that do not contribute
to the polar locus of the local zeta function (Theorem
\ref{unmask}). We prove a result concerning the limit of the zeta
function at $-\infty$ in Section \ref{degree}, Theorem
\ref{theoremdegree}. To conclude, in Section \ref{proofs}, we
prove the generalized Monodromy and Holomorphy Conjecture in the
case where $dim(X)=2$.
\section{Several variables Igusa local zeta functions}\label{def}

\subsection{}Let us fix some notation.
We choose a prime number $p$.
We denote by $K$ a finite extension of $\QQ_p$,
with ring of integers $R$, and residue field $k$, a finite extension
of $\FF_p$. We choose
an uniformizing parameter $\varpi$ and denote by $q$ the cardinality of $k$.
We denote by ${\rm val}$ the valuation on $K$ and
we normalize the norm on
$K$ by setting $|x| = q^{-{\rm val} (x)}$.
For $x$ non zero in $K$, we define $\ac (x)$ as the image in
$k$ of $x \varpi^{-{\rm val} (x)}$ and we set
$\ac (x) = 0$ if $x = 0$.

Let $X$ be a smooth compact $p$-adic (locally) analytic
variety (hence $X$ is a finite disjoint union of balls) of pure dimension
$d$ and let $\omega$ be a gauge form on
$X$, i.e. a nowhere vanishing analytic $n$-form on
$X$. In a standard way (cf. \cite{serre}),
we may attach to $\omega$ a volume form
$|\omega|$ on $X$.

\subsection{}We denote by $\cC (\GG_m^r (k), \CC)$
the group
of characters
$$\chi = (\chi_1, \dots, \chi_r) :
\GG_m^r (k) = (k^*)^r \longrightarrow \CC^*.$$
For any analytic function
$F = (f_1, \dots, f_r) : X \rightarrow K^r$,
for
$\phi$ locally constant on $X$, for $\chi$ in
$\cC (\GG_m^r (k), \CC)$
and for $s = (s_1,\dots, s_r) $ in
$\CC^r$,
we consider Igusa's local zeta function
$$
Z_{\phi} (F, \chi ; s)
:=
\int_X \phi \prod_{i = 1}^r
\chi_i (\ac (f_i (x)))
\prod_{i = 1}^r
|f_i (x)|^{s_i} |\omega|,
$$
which is defined when $\Re (s_i) > 0$, for each $i$, where $\Re
(s_i)$ denotes the real part of $s_i$. We will often omit
$\phi,F,\chi$ from the notation, and denote the zeta function
simply by $Z(s)$. Furthermore, when $X=R^{n}$, and $\phi$ is the
characteristic function of $P^{n}$, with $P$ the maximal ideal in
$R$, we write $Z_0(s)$ instead of $Z_{\phi}(F,\chi;s)$.

\subsection{}Consider
an embedded resolution $\pi : Y \rightarrow X$ of the $f_j$'s. By
this we mean that $\pi : Y \rightarrow X$ is a proper morphism
between compact locally analytic $p$-adic varieties, that $\pi$
induces an isomorphism between $Y \setminus \pi^{-1} (D)$ and $X
\setminus D$, with $D = \cup_{j} f_j^{-1} (0)$, and that $\pi^{-1}
(D)$ is a divisor $E$ with (global) normal crossings. We denote by
$E_i$, $i \in I$, the set of irreducible components of $E$. In
particular, the $E_i$ are smooth. We denote by $N_i^j$ the order
of vanishing of $f_j \circ \pi$ at the generic point of $E_i$, and
by $\nu_i - 1$ the order of vanishing of $\pi^* \omega$ at the
generic point of $E_i$. We shall write $N_i s$ for $\sum_j N_i^j
s_j$, and $\chi^{N_i}$ for $\prod_j \chi_j^{N_i^j}$.

\begin{theorem}[Igusa]\label{rat}
The function
$Z_{\phi} (F, \chi ; s)$ is a rational function
in the variables $q^{-s_j}$. Furthermore,
given
an embedded resolution $\pi : Y \rightarrow X$  of the $f_j$'s,
a necessary condition for
$Z_{\phi} (F, \chi ; s)$ to have a pole for a given value
of $\chi$ and $s$ is that
$$
\prod_{i \in I} (1 -q^{- \nu_i} q^{ -N_i s}) = 0
$$
and
$$
\prod_{i \in I} (\chi^{N_i} - 1) =0.
$$
\end{theorem}

\subsection{}By Theorem \ref{rat},
the set of real parts of the poles of $Z_{\phi} (F, \chi ; s)$ are
located on a finite set of affine hyperplanes
$$
H : \quad \sum_j N^j s_j + \nu = 0
$$
with $N^j$ in $\NN$ and $\nu > 0$. We denote this finite set of
hyperplanes by $\cP (Z_{\phi} (F, \chi ; s))$.

\section{Alexander modules}\label{alex}

\subsection{Notations and conventions}
We fix an algebraically closed field $k$ of characteristic
exponent $p$ (i.e. $p = 1$ when $k$ is of characteristic $0$, and
$p$ is equal to the characteristic of $k$ when $k$ is of
characteristic $> 0$). We fix a prime number $l \not= p$. We
denote by $\bar \QQ_{l}$ a fixed algebraic closure of $\QQ_{l}$.
We set $\hat \ZZ (1) (k) := \limproj_{(p, n) = 1} \mu_n (k)$ and
$\ZZ_{l} (1) (k) := \limproj_{n} \mu_{l^n} (k)$. Let $T$ be a
$k$-torus, $T \simeq \GG_{m, k}^r$. We set $X_{\ast} (T) := {\rm
Hom}_{k-{\rm gr}} (\GG_{m, k}, T)$, $\hat X_{\ast} (T) := X_{\ast}
(T) \otimes_{\ZZ} \hat \ZZ (1) (k)$ and $X_{\ast l} (T):= X_{\ast}
(T)\otimes_{\ZZ} \ZZ_{l} (1) (k)$. We denote by $\pi_1 (T)$ the
fundamental group of $T$ pointed at 1, by $\pi_1 (T)^t$ the tame
fundamental group of $T$, and by $\pi_1 (T)_{l}$ the maximal
pro-$l$ quotient of $\pi_1 (T)^t$. By Kummer theory, there are
canonical isomorphisms $\hat X_{\ast} (T) \simeq \pi_1 (T)^t$ and
$X_{\ast l} (T) \simeq \pi_1 (T)_{l}$.

\subsection{Torsion by the generic character}
Let $R$ be a complete regular local ring with residual
characteristic $l$. For every $k$-torus $T$, we set
$$
\Omega^R_T := R [[\pi_1 (T)_{l}]].
$$
This is a complete regular local ring of residual characteristic
$l$. By Kummer theory, $\Omega^R_T$ may be identified with $R
[[X_{\ast l} (T)]]$. For every choice of a generator $\gamma$ of
$\ZZ_{l} (1) (k)$, we have, by Iwasawa theory, an isomorphism
$$
R [[\ZZ_{l} (1) (k)]] \simeq R[[t]]
$$
sending $\gamma$ to $1 +t$.
Hence, for every isomorphism $T \simeq (\GG_{m, k})^r$,
we have, $\gamma$ being given,
an isomorphism
$$
\Omega^R_T \simeq R[[t_1, \ldots, t_r]].
$$
There is a canonical continuous character
$$
{\rm Can}_T : \pi_1 (T) \rightarrow ({\Omega^R_T})^{\times}
$$
obtained by composition of the tautological character $\pi_1
(T)_{l} \rightarrow ({\Omega^R_T})^{\times}$, sending $\gamma$ to
$\gamma$, with the projection $\pi_1 (T) \rightarrow \pi_1
(T)_{l}$. We shall denote by $L_T^R$ the corresponding
$\Omega_T^R$-sheaf on $T$. It is a tamely ramified lisse
constructible free twisted rank 1 $\Omega_T^R$-sheaf on $T$.

\subsection{The scheme $\cC (T)$}
 In \cite{GabLoe}, Gabber and Loeser defined a $\bar \QQ_{l}$-scheme $\cC (T)$
whose set of closed points may be canonically identified with the
group of continuous characters $\pi_1 (T)^t \rightarrow \bar
\QQ_{l}^{\times}$. It is the disjoint union of the translates of
the connected component $\cC (T)_{l}$ of the trivial character by
the group of characters $\chi : \pi_1 (T)^t \rightarrow \bar
\QQ_{l}^{\times}$ of finite order prime to $l$. By its very
definition,
$$
\cC (T)_{l} := \Spec \bar \QQ_{l} \otimes \ZZ_{l} [[\pi_1
(T)_{l}]].
$$
By \cite{GabLoe} 3.2.2, $\cC (T)_{l}$ is a noetherian regular
scheme.

\subsection{}\label{cot}
If $\pi : T \rightarrow T'$
is a morphism of tori, we denote
by $\pi^{\vee} : \cC (T') \rightarrow \cC (T)$
the morphism which one deduces by
functoriality.
We say that
$\pi : T \rightarrow T'$
is a quotient of $T$, if $\pi$ is surjective with kernel a torus.
We call a subscheme of
$\cC (T)$
of the form $\pi^{\vee} (\cC (T'))$, with $\pi : T \rightarrow T'$
a quotient of $T$,
a cotorus of $\cC (T)$. A subscheme of the form
$ \chi \cdot \pi^{\vee} (\cC (T'))$,
with $\chi$ a closed point of $\cC (T)$, is called a
translated cotorus.
If $Z$ is cotorus translated by a character of order $n$
(i.e. $Z = \chi \cdot \pi^{\vee} (\cC (T'))$
with $\chi$ of order $n$),
we
set
$$
\overline Z := \bigcup_{1 \leq i \leq n} \chi^i \cdot \pi^{\vee} (\cC (T')).
$$
More generally, if $Z$ is a finite union $\cup Z_i$ of cotori
translated by a character of finite order, we set $\overline Z :=
\cup_i \overline Z_i$.

\subsection{}\label{fo}Assume $k$ is the
algebraic closure of a finite field $k_0$ with $q$ elements. Fix
an isomorphism $\bar \QQ_{l} \simeq \CC$. Then we may identify
$\cC (\GG_m^r (k_0), \CC)$ with the set of closed points of $\cC
(\GG_{m, k}^r)$ of order dividing $q - 1$, since we have a
canonical morphism $\hat \ZZ (1) (k) \rightarrow \mu_{q - 1} (k)
\simeq k_0^{\ast}$.

\subsection{}\label{rest}Assume $K$ is the algebraic closure of some local field
having $k$ as residue field.
The canonical injection $\hat \ZZ (1) (k) \hookrightarrow
\hat \ZZ (1) (K)$
induces an epimorphism
$r : \cC (\GG_{m, K}^r) \rightarrow \cC (\GG_{m, k}^r)$.

\subsection{Alexander modules}
Now we consider the affine space $$\A^r_k := \Spec k [T_1, \ldots,
T_r]$$ and the torus $j : T = \Spec k [T_1, T_1^{-1},\ldots, T_r,
T_r^{-1}] \hookrightarrow \A^r_k$. Let $F : X \rightarrow \A^r_k$
be a morphism of $k$-varieties and set $j : X_T := F^{-1} (T)
\hookrightarrow X$.

Let $R$ be the ring of integers of a finite extension of
$\QQ_{l}$, and let $\cF$ be an object in $D^b_c (X_T, R)$. As
explained in  \cite{GabLoe} A.1, one can define in a natural way
objects $\cF \otimes_R F^{\ast} L_T^R$ and $Rj_{\ast}  (\cF
\otimes_R F^{\ast} L_T^R)$ in $D^b_c (X_T, \Omega^R_T)$ and $D^b_c
(X, \Omega^R_T)$,  respectively. Hence, for every geometric point
$\xi$ of $X \setminus X_T$, and every integer $q \geq 0$, one may
consider
$$
\cA^{q, R}_{F, \xi} (\cF) :=
R^qj_{\ast}  (\cF \otimes_R F^{\ast} L_T^R)_{\xi}
$$
and view it as a coherent $\Omega^R_T$-module.

Now, if $\cF$ is an object in $D^b_c (X_T, \bar \QQ_{l})$, we
define a coherent $\cC (T)$-module $\cA^{q}_{F, \xi} (\cF)$, the
$q$-th Alexander module  of $\cF$ at $\xi$, as follows. We choose
a model $\cF_R$ in $D^b_c (X_T, R)$ for some $R$, and ask for  the
restriction of $\cA^{q}_{F, \xi} (\cF)$ to the connected component
$\cC (T)_{l}$ to be equal to
$$\cA^{q}_{F, \xi} (\cF)_{|\cC (T)_{l}}
:= \cA^{q, R}_{F, \xi} (\cF_R) \otimes \cC (T)_{l},$$ which is
clearly independent of the chosen model $\cF_R$. For  the
restriction of $\cA^{q}_{F, \xi} (\cF)$ to a general connected
component of the form $\{\chi\} \times \cC (T)_{l}$, with $\chi$ a
character of finite order prime to $l$, we set
$$\cA^{q}_{F, \xi} (\cF)_{|\{\chi\} \times \cC (T)_{l}}
:= \cA^{q}_{F, \xi} (\cF \otimes \cL_{\chi})_{|\cC (T)_{l}},
$$
with $\cL_{\chi}$ the Kummer sheaf attached to $\chi$ (its
definition is similar to the construction of
$\mathcal{L}_{\psi,\alpha}$ in Section \ref{form}).

If $\cF$ belongs to $D^b_c (X, \bar \QQ_{l})$, we shall write
$\cA^{q}_{F, \xi} (\cF)$ for $\cA^{q}_{F, \xi} (j^{\ast}\cF)$.
Also, if $k$ is not algebraically closed, we shall still write
$\cA^{q}_{F, \xi} (\cF)$ for the Alexander module attached to $F
\otimes \bar k  : X \otimes \bar k \rightarrow \A^n_{\bar{k}}$.

\subsection{Support of Alexander modules}

\begin{prop}\label{support}
Let $F : X \rightarrow \A^n_k$ be a morphism of $k$-varieties, and
let $\cF$ be an object in $D^b_c (X_T, \bar \QQ_{l})$. Then, for
every geometric point $\xi$ of $X \setminus X_T$ and every integer
$q$, the support of the module $\cA^{q}_{F, \xi} (\cF)$ is
contained in a finite set, not depending on $\xi$ and $q$, of
translated cotori of codimension 1 of $\cC (T)$. Furthermore, when
$\cF = \bar \QQ_{l}$, the support of the module $\cA^{q}_{F, \xi}
(\bar \QQ_{l})$ is contained in a finite set, not depending on
$\xi$ and $q$, of cotori of codimension 1 of $\cC (T)$ translated
by characters of finite order.
\end{prop}
\begin{proof}
We may suppose that $\cF$ is the extension by zero of a lisse
sheaf on a smooth irreducible locally closed subvariety $U$ of
$X_T$. We may furthermore assume that $U$ is dense in $X_T$, and
that the monodromy of $\cF$ is pro-$l$; see \cite{GabLoe}, Prop.
4.3.1'.

By de Jong's theorem \cite{dJ}, there exists a generically \'etale
alteration $\pi:X'\rightarrow X$, where $X'$ is regular, and the
inverse image $\pi^{-1}(Z)$ is a normal crossing divisor, with
$Z=X-U$. Shrinking $X$ if necessary, we may suppose that $\pi$ is
\'etale over $U$.

\begin{picture}(100,100)
\put(120,70){$\pi^{-1}(U)$} \put(160,74){\vector(1,0){40}}
\put(170,76){$u'$} \put(205,70){$\pi^{-1}(X_T)$}
\put(250,74){\vector(1,0){40}} \put(260,77){$j'$}
\put(295,70){$X'$}

\put(135,65){\vector(0,-1){40}} \put(220,65){\vector(0,-1){40}}
\put(300,65){\vector(0,-1){40}} \put(303,45){$\pi$}

\put(130,10){$U$} \put(150,14){\vector(1,0){50}}\put(170,16){$u$}
\put(215,10){$X_T$}
\put(240,14){\vector(1,0){50}}\put(260,17){$j$} \put(295,10){$X$}
\end{picture}

For any object $\mathcal{G}$ in $D^{b}_{c}(X_T,R)$, and any
character $\chi$ of finite order prime to $l$, we have a short
exact sequence
$$0\rightarrow u_{!}u^{*}\mathcal{G}\otimes \mathcal{L}_{\chi}\rightarrow \mathcal{G}\otimes \mathcal{L}_{\chi}\rightarrow i_{*}i^{*}\mathcal{G}\otimes \mathcal{L}_{\chi}\rightarrow 0 \,,$$
where $i$ is the inclusion of $Z\cap X_T$ in $X_T$, and
$\mathcal{L}_{\chi}$ is the Kummer sheaf of $\chi$ on $X_T$.
Applying the derived functor $R^{q}j_{*}$ yields a long exact
sequence, hence, it suffices to investigate the support of the
Alexander complex of $u_{!}u^{*}\mathcal{G}$ and
$i_{*}i^{*}\mathcal{G}$.

The theorem holds for the latter, by Noetherian induction, since
$i_{*}$ is exact, hence, by Leray's spectral sequence,
$$R^{q}j_{*}(i_{*}i^{*}\mathcal{G}\otimes \mathcal{L}_{\chi})\cong R^{q}(j\circ
i)_{*}(i^{*}\mathcal{G}\otimes i^{*}\mathcal{L}_{\chi})\,.$$

As for the former, by proper base change, one sees that
$$u_{!}u^{*}\mathcal{G}\otimes \mathcal{L}_{\chi}\cong \pi_{*}u'_{!}\pi^{*}u^{*}\mathcal{G}\otimes \mathcal{L}_{\chi} \,.$$
Leray's spectral sequence reads
$$E_{2}^{p,q}= (R^{p}j_{*})(R^{q}\pi_{*})(u'_{!}\pi^{*}u^{*}\mathcal{G}\otimes \pi^{*}\mathcal{L}_{\chi})
\Rightarrow R^{p+q}(j\circ
\pi)_{*}(u'_{!}\pi^{*}u^{*}\mathcal{G}\otimes
\pi^{*}\mathcal{L}_{\chi})\,.$$ We use the proper base change
theorem to compute the stalks of
$(R^{q}\pi_{*})(u'_{!}\pi^{*}u^{*}\mathcal{G}\otimes
\pi^{*}\mathcal{L}_{\chi})$. These are zero outside of $U$, and
since $\pi$ is \'etale over $U$, they vanish in all other points
as well, if $q>0$. Therefore, our spectral sequence yields
isomorphisms
$$R^{p}j_{*}(u_{!}u^{*}\mathcal{G}\otimes \mathcal{L}_{\chi})\cong R^{p}(\pi\circ j')_{*}(u'_{!}\pi^{*}u^{*}\mathcal{G}\otimes \pi^{*}\mathcal{L}_{\chi})\,.$$
After applying Leray's spectral sequence and proper base change
once more, one sees that it suffices to consider the support of
$R^{p}j'_{*}(u'_{!}\pi^{*}u^{*}\mathcal{G}\otimes
\pi^{*}\mathcal{L}_{\chi})_{\xi'}$, where $\xi'$ is a geometric
point of $\pi^{-1}(X\setminus X_T)$.

This means we may assume that $W$ is smooth over $k$, that the
complement of $U$ in $X$ is a normal crossing divisor, and that
the monodromy of $\mathcal{F}$ on $U$ is pro-$l$. Now one can
proceed as in the proof of \cite{GabLoe}, Proposition 4.3.1'.
\end{proof}

\section{Conjectures}\label{conj}
\subsection{}In this section, we shall consider
a number field $K$, and a smooth algebraic variety $X$ of
dimension $d$ over $K$, together with a gauge form $\omega$ on
$X$. Let $F := (f_1, \dots, f_r)$ be a morphism $X \rightarrow
\A^r_K$. For every finite place $\gP$ of $K$, we denote by
$K_{\gP}$ the completion of $K$ and by $R_{\gP}$ its ring of
integers. Taking an appropriate model $\cX$ of $X$ over the
integers of $K$ imposes, at least for almost all $\gP$, the
structure of a smooth compact locally analytic $K_{\gP}$-variety
on $X_{\gP} := \cX (R_{\gP})$. Also, taking appropriate models of
$\omega$ and $F$, we obtain a gauge form $\omega_{\gP}$ on
$X_{\gP}$ and an analytic function $F_{\gP} = (f_{1 \gP}, \dots,
f_{r \gP}) : X_{\gP} \rightarrow K_{\gP}^r$. All these data are
independent of choices, for almost all $\gP$. Hence, for suitable
$\chi$ and $\phi$, we can consider Igusa's zeta function $Z_{\phi}
(F_{\gP}, \chi ; s)$.

\subsection{}Let us first review
the $r = 1$ case, i.e. the case when there is only one function $f$.
We write $f$ for $F = (f)$.

We shall say a locally constant function $\phi$ on $X_{\gP} := \cX
(R_{\gP})$ is residual if its value at a point $x$ depends only of
the image in $\cX (\FF_{\gP})$, where $\FF_{\gP}$ denotes the
residue field.

\begin{monoconjecture}[Igusa, \cite{DenefBour}]\label{monoconjecture}
For almost all finite places $\gP$, if $s_0$ is a pole of
$Z_{\phi} (f_{\gP}, \chi ; s)$, with $\phi$ residual on $X_{\gP}$,
then $\exp (2 \pi i {\rm Re} (s_0))$ is an eigenvalue of the local
monodromy at some point of $f^{-1} (0)$.
\end{monoconjecture}

\begin{holoconjecture}[Denef, \cite{Denef5}]\label{holoconjecture}
For almost all finite places $\gP$, if $Z_{\phi} (f_{\gP}, \chi_0
; s)$, with $\phi$ residual on $X_{\gP}$, is not a holomorphic
function of $s$, then the order of $\chi_0$ divides the order of
some eigenvalue of the local monodromy at some point of $f^{-1}
(0)$.
\end{holoconjecture}

Here we should make more precise what is meant by ``some
eigenvalue of the local monodromy at some point of $f^{-1} (0)$''.
The easiest way is to fix a complex embedding of $K$ in $\CC$ and
to stick to  the complex points of $X$ with complex topology. We
can consider the usual complex of nearby cycles $R \psi \CC$ on
$f^{-1} (0)$, and then ``some eigenvalue of the local monodromy at
some point of $f^{-1} (0)$'' means ``some eigenvalue of the
monodromy action on  some cohomology object of the stalk $(R \psi
\CC)_x$ at some complex point $x$ of $f^{-1} (0)$''. More
algebraically, and this amounts to the same by the classical
comparison theorems, one could consider the complex of nearby
cycles $R \psi \QQ_{l}$ on $f^{-1} (0)$ for the \'etale topology
and $l$-adic cohomology, and
 interprete it as
``some eigenvalue of the geometric monodromy action on  some
cohomology object of the stalk $(R \psi \QQ_{l})_x$ at some
geometric point $x$ of $f^{-1} (0)$''.

\subsection{}We fix an isomorphism $\bar \QQ_{l} \simeq \CC$.
In particular, we have an embedding of $\bar \QQ$ in $\CC$. We
denote by $\gamma$ the topological generator $\gamma = (\zeta_n)$,
$\zeta_n = \exp (2 \pi i / n)$, of $\hat \ZZ (\QQ)$. To an affine
hyperplane $H$ with equation $\sum_j N^j s_j + \nu = 0$ in
$\CC^r$, with $N^j$ and $\nu$ in $\NN$, $\nu > 0$, we associate
$\cH$, a finite union of cotori translated by a character of
finite order in $\CC (\GG_{m, \bar K})$, as follows. We deduce
from $\gamma$ a topological basis $(\gamma_1, \dots, \gamma_r)$ of
$\pi_1 (\GG_{m, \bar K}^r)^t = \hat \ZZ (1)^r (\bar K)$. For any
$\alpha$ in $(\QQ/ \ZZ)^r$, we denote by $\varphi (\alpha)$ the
character which sends $(\gamma_1, \dots, \gamma_r)$ to $\exp (2
\pi i \alpha)$. We consider the image $[H (\QQ)]$ of $H (\QQ)$ in
$(\QQ/ \ZZ)^r$, and we define $\cH$ as the Zariski closure of
$\{\varphi (\alpha)|\alpha \in [H (\QQ)]\}$ in $\CC (\GG_{m, \bar
K}^r)$. It is clearly a finite union of cotori translated by a
character of finite order.

\subsection{}Fran\c cois Loeser proposes the following generalisation of the monodromy
conjecture:

\begin{gmonoconjecture}\label{gmonoconjecture}\
For almost all finite place $\gP$, if $H$ is a polar hyperplane in
$\cP (Z_{\phi} (F_{\gP}, \chi ; s))$, with $\phi$ residual on
$X_{\gP}$, then, for every irreducible component $Z$ of $\cH$,
there exists an integer $q$ and a geometric point $\xi$ such that
$Z$ is contained in the support of $\cA^{q}_{F,  \xi} (\bar
\QQ_{l})$.
\end{gmonoconjecture}

Loeser also generalizes the holomorphy conjecture as follows
(using notations and conventions from \ref{fo} and \ref{rest}):

\begin{gholoconjecture}\label{gholoconjecture}
For almost all finite place $\gP$, if $Z_{\phi} (F_{\gP}, \chi_0 ;
s)$ with $\phi$ residual on $X_{\gP}$ is not a holomorphic
function of $s$, then there exists an integer $q$ and a geometric
point $\xi$, such that $\chi_0$ belongs to $r (\overline{{\rm
Supp} \, \cA^{q}_{F,  \xi} (\bar \QQ_{l})})$.
\end{gholoconjecture}

\begin{remark}By Proposition \ref{monoco}, when $r = 1$,
Conjectures \ref{gmonoconjecture} and \ref{gholoconjecture}
are indeed equivalent to
Conjectures \ref{monoconjecture} and \ref{holoconjecture}.
\end{remark}

\section{Comparison with the transcendent theory}\label{trans}
In this section, we compare our definition of the Alexander
complex with the classical definition \cite{Sabbah} in the complex
case, following Deligne's treatment of the nearby cycles
\cite{SGA7b}. Let $X$ be a complex manifold, and let
$F=(f_1,\ldots,f_{r}):X\rightarrow \C^{r}$ be an analytic map. We
denote the complex analytic torus $\mathbb{G}_{m,\C}^{r}(\C)$ by
$T$, the inverse image of $T$ under $F$ by $X_{T}$, and the
complement of $X_{T}$ in $X$ by $X_0$. Let
$\pi:\tilde{T}\rightarrow T$ be a universal cover of $T$. Consider
the following diagram:

\begin{picture}(200,80)(-30,20)
\put(30,80){$X_0$} \put(45,83){\vector(1,0){30}} \put(82,80){$X$}
\put(125,83){\vector(-1,0){30}} \put(130,80){$X_T$}
\put(175,83){\vector(-1,0){28}} \put(180,80){$X_{T}\times_T
\tilde{T}=X_{\tilde{T}}$}

 \put(86,77){\vector(0,-1){30}}
\put(137,77){\vector(0,-1){30}} \put(196,77){\vector(0,-1){30}}

\put(90,58){$F$} \put(162,42){$\pi$} \put(162,85){$\tilde{\pi}$}
\put(57,85){$i$} \put(107,86){$j$}

\put(82,35){$\C^{r}$} \put(127,40){\vector(-1,0){32}}
\put(135,35){$T$} \put(185,40){\vector(-1,0){40}} \put(193,35){$
\tilde{T}$}

\end{picture}

\noindent In the univariate case ($r=1$), after reducing $X$ to a
sufficiently small neighbourhood of a point $x\in F^{-1}(0)$ (i.e.
the intersection of a small ball centered at $x$ with the inverse
image of a small ball around the origin in $\C$), the space
$X_{T}\times_T \tilde{T}$ is called the canonical Milnor fiber of
$F$ at $x$, since it is homotopy-equivalent to the fiber of the
Milnor fibration at $x$ \cite{Milnor}\cite{Kuli}. However, in the
multivariate case, $F$ will, in general, not be a fibration in any
neighbourhood of $x$.

Let $\mathcal{F}$ be an object in $D_{c}^{b}(X_T,\C)$. We define
the Alexander complex to be
$$R\psi_{F}(\mathcal{F})=i^{*}R(j\circ \tilde{\pi})_{*}
\tilde{\pi}^{*}\mathcal{F}\,.$$ The deck transformations of
$X\times_{T}\tilde{T}$ make $R\psi_{F}(\mathcal{F})$ into a
constructible $\C[\Z^{r}]$-module.

Now let $A$ be a regular Noetherian ring, $\mathcal{F}$ an
$A$-constructible sheaf on $X_0$, and $\mathcal{G}$ an object in
$D_{c}^{b}(X_0,A)$. The support Supp$_x \mathcal{F}$ of
$\mathcal{F}$ at a point $x$ of $X_0$ is defined to be the
collection of prime ideals $\mathcal{P}$ of $A$ such that
$\mathcal{F}_{x}\otimes_{A}A_{\mathcal{P}}\neq 0$; the support of
$\mathcal{G}$ at $x$ is the union of Supp$_{x}h^{i}(\mathcal{G})$.
If $d$ is the dimension of the support of $\mathcal{G}$ at $x$, we
denote by Supp$_{x,d}\,\mathcal{G}$ the collection of prime ideals
in Supp$_{x}\mathcal{G}$ of coheight $d$. The cycle associated to
$\mathcal{G}$ at $x$ is the cycle
$$\zeta_{x,d}(\mathcal{G})=\sum_{\mathcal{P}\in\,
\mathrm{Supp}_{x,d}\mathcal{G}}\chi_x(\mathcal{G}_{\mathcal{P}})V(\mathcal{P}),$$
where $V(\mathcal{P})$ is the subvariety of Spec\,$A$ associated
to $\mathcal{P}$, $\mathcal{G}_{\mathcal{P}}$ is the localization
$\mathcal{G} \otimes_A A_{\mathcal{P}}$, and
$\chi_x(\mathcal{G}_{\mathcal{P}})$ equals $\sum
(-1)^{i}\mathrm{length}\,h^{i}(\mathcal{G_P})_x$. In our setting,
$A$ equals $\C[t_i,t_i^{-1}]_{i=1}^{r}$, and $d$ will be equal to
$r-1$, so we can replace the cycle by a defining function. This
zeta function $\zeta_{x,d}(\mathcal{G})$ is a rational function in
$t_1,\ldots,t_r$, defined modulo multiplication with a monomial.
We get a duality formula
$$\zeta_{x,d}(D_A\mathcal{G})(t_1,\ldots,t_r)=\zeta_{x,d}(\mathcal{G})(t_1,\ldots,t_r)^{-1},$$
where $D_A\mathcal{G}$ is the dual complex of $\mathcal{G}$.

In his article \cite{Sabbah}, Sabbah formulates a generalized
A'Campo formula \cite{A'C}: let
$F=(f_1,\ldots,f_r):(\C^{n},0)\rightarrow (\C^{r},0)$ be an
analytic germ, and let $h:E\C^{n}\rightarrow \C^{n}$ be a proper
modification, where $E\C^{n}$ is smooth. Assume that there exists
a divisor $E=\cup E_i$ on $E\C^{n}$ with normal crossings, such
that $(f_i\circ h)^{-1}(0)$, as well as $h^{-1}(0)$ are unions of
components of $E$. Then
$$\zeta_0(R\psi_F\C)(t_1,\ldots,t_r)=\prod_{\{i\in I\,|\,E_i\in
\pi^{-1}(0)\}}(t^{N_i}-1)^{-\chi(E^{o}_i)}\,,$$ where
$t^{N}=t_1^{N^{1}}\ldots t_r^{N^r}$, the $N_i$ are multiplicity
vectors of $F$ along the $E_i$, and $E_i^{o}$ is equal to the set
of points on $E_i$ which do not belong to any other divisor $E_j$.

It is well-known \cite{Sabbah} that $Rh_*\circ R\psi_{F\circ
h}=R\psi_{F}\circ Rh_*$. The comparison of the analytic with the
algebraic definition is based upon an explicit computation of the
Alexander complex associated to $F\circ h$, which reduces to the
computation of the Alexander complex associated to
$$G:Y=\C^{n}\rightarrow \C^{r}:x=(x_1,\ldots,x_n)\mapsto
(x^{N^{1}},\ldots,x^{N^{r}})\ ,$$ with $N^{j}\in \N^{n}$. We may
assume that for each index $j=1\ldots n$, there exists a $k$ such
that $N^{k}_j\neq 0$. The map
$$p_1:\C^{n}\rightarrow Y_T=(\C_0)^{n}:(w_1,\ldots,w_n)\mapsto
(e^{2i\pi w_1},\ldots, e^{2i\pi w_n})$$ is a universal covering
space. The space $Y_{\tilde{T}}$ falls apart in $\kappa$ connected
components, where $\kappa$ equals $\prod_i
gcd\,\{N^{j}_{i}\}_{j=1}^{r}$. These components are transitively
permutated by the monodromy actions. For each of these components
$Z$, we get an induced covering projection $p_2:\C^{n}\rightarrow
Z$ with
$$M=\{(z_1,\ldots,z_n)\in \Z^{n}\,|\,\sum z_i \,N^{j}_i=0\
\mathrm{for\  each}\ j\}\subset \Z^{n}=\pi_1(Y_T)$$ as
transformation group. This means that, for each constant sheaf
$\Lambda$ on $Y$ and each $q$,
$$R^{q}\psi_{G}(\Lambda)_0=H^{q}(Y_{\tilde{T}},\Lambda)=\wedge^{q}(M\otimes \Lambda)\otimes \psi^{0}_{G}(\Lambda)_0\,,$$ which
follows from the fact that, for each sheaf $\mathcal{F}$ on $Z$,
the global sections of $\mathcal{F}$ are the global sections of
$p_2^{*}\mathcal{F}$ which are invariant under $M$, and the
resulting spectral sequence.

Let us try to reformulate this definitions in a more algebraic
way. For each $i=1,\ldots,r$, let $S_i$ be the spectrum of a
strictly Henselian ring with residue field $k$, closed point $s_i$
and generic point $\eta_i$. Consider a geometric point
$\bar{\eta_i}$, localized at $\eta_i$. The residue field of the
closed point of the normalization $\bar{S}_i$ of $S_i$ in
$k(\bar{\eta_i})$ is an algebraic separably closed extension
$k(\bar{s}_i)$ of $k$, and thus defines a geometric point
$\bar{s}_i$ localized at $s_i$. Let $F:X\rightarrow
S$ be a scheme of finite type over the fiber
product $S$ of the $S_i$ over Spec\,$k$. Fix a subset $J$ of
$\{1,\ldots,r\}$. Let $\Lambda$ a Noetherian torsion ring in which
the characteristic exponent of $k$ is invertible, and let
$\mathcal{F}$ be an object in $D_{c}^{b}(X_T,\Lambda)$.

Consider the diagram

\begin{picture}(200,80)(-30,20)
\put(69,80){$X_{s_J}$} \put(85,83){\vector(1,0){40}}
\put(130,80){$X$} \put(185,83){\vector(-1,0){42}}
\put(190,80){$X_{\bar{\eta}}$}

 \put(76,76){\vector(0,-1){30}}
\put(134,77){\vector(0,-1){31}} \put(198,76){\vector(0,-1){30}}

\put(139,58){$F$} \put(162,85){$\tilde{\pi}$}
 \put(107,85){$i$}

\put(73,35){$s_{J}$}
 \put(85,38){\vector(1,0){42}} \put(131,35){$S$}
\put(193,38){\vector(-1,0){20}} \put(197,35){$ \bar{\eta}$}
\put(163,35){$\eta$} \put(159,38){\vector(-1,0){20}}

\end{picture}

\noindent where $\eta=(\eta_i)$, and $s_J=(\alpha_i)$ with
$\alpha_i=s_i$ if $i\in J$ and $\alpha_i=\eta_i$ else. We define
the $J$-th Alexander complex to be
$$R\psi_{F,J}(\mathcal{F})=i^{*}R\tilde{\pi}_{*}\mathcal{F}_{\bar{\eta}}\,.$$
It is endowed with an action of $\prod
\mathrm{Gal}(\bar{\eta_i}/\eta_i)$; the $i$-part acts trivially if
$i\notin J$. If we replace, for each $i$, $\bar{\eta_i}$ by the
point $\eta_{i,t}$ associated to the maximal tamely ramified
extension of $k(\eta_i)$ in $k(\bar{\eta_i})$, we get the tame
Alexander complex $R\psi_{F,J}^{(t)}$.

 It follows from proper base change that $Rh_*\circ R\psi_{F\circ h,J}^{(t)}=R\psi_{F,J}^{(t)}\circ
 Rh_*$ if $h:Y\rightarrow X$ is a proper morphism. In our
 comparison theorem, we will apply this property to the case where
 $h$ is a simultaneous embedded resolution for the $f_i$; so let
 us investigate what happens if, for each $i\in J$, $X_{\eta_i}$
 is smooth over $\eta_i$, and $D=\cup_{i\in J}X_{s_i}$ is a normal
 crossing divisor. The calculation follows the same lines as in
 the analytic case; see \cite{SGA7a} for the univariate case $r=1$.

  Let $\bar{x}$ be a
 geometric point of $X$, localized at a point $x$ mapped to $s_J$
 by $F$. Consider the strict localization $X_{(\bar{x})}$ of $X$
 at $\bar{x}$, let $U$ be equal to the complement of the inverse image of $D$ in
 $X_{(\bar{x})}$, and let $t_j$ be a local equation for $D_j$. Now
 put $X_n=X_{(\bar{x})}[t_j^{\frac{1}{n}}]$, and denote by $U_n$ the
 inverse image of $U$ in $X_n$. By purity and the long exact
 sequence for relative cohomology,
 $H^{q}(U_n,\Lambda)=\wedge^{q}\Lambda(-1)^{C}$, where $C$ is the
 number of components $D_j$ containing $x$. For $m=nd$, the
 morphism from $H^{q}(U_n,\Lambda)$ to $H^{q}(U_m,\Lambda)$,
 associated to the obvious morphism $U_m\rightarrow U_n$, is
 induced by multiplication by $d^{q}$. If we define $\tilde{U}$ to
 be the projective limit of the $U_n$, taken over the indices $n$
 prime to the characteristic exponent of $k$, we get
 $$\left\{\begin{array}{ll}H^{0}(\tilde{U},\Lambda)=\Lambda &
 \\ H^{q}(\tilde{U},\Lambda)=\lim_{\to}H^{q}(U_n,\Lambda)=0 &
 \mbox{if $q\neq 0$.}\end{array}\right.$$
 The scheme $\tilde{U}$
 is a procover of $U$ with group $\hat{\Z}(1)(k)^{C}$.
 Furthermore, it is a procover with group $M\otimes
 \hat{\Z}(1)(k)$ of each of the connected components of
 $U_{\eta_t}$, which are permutated transitively by the tame
 inertia group $I_t=\hat{\Z}(1)(k)^r$. Since, by definition,
 $R\psi^{(t),q}_{F,J}(\Lambda)_{\bar{x}}=H^{q}(U_{\eta_t},\Lambda)$,
 the Hochschild-Serre spectral sequence yields
 $$R^{q}\psi^{(t)}_{F,J}(\Lambda)_{\bar{x}}\cong \wedge^{q}(M\otimes
 \Lambda(-1))\otimes_{\Lambda}R^{0}\psi^{(t)}_{F,J}(\Lambda)_{\bar{x}}\,,$$
 and $R^{0}\psi^{(t)}_{F,J}(\Lambda)_{\bar{x}}$ is a
 $\Lambda$-algebra $\Lambda^{E}$, with $E$ a set on which the tame
 inertia group acts transitively, and $|E|$ equal to the greatest
 divisor of $\kappa$ prime to $p$. Of course, $M$ and $\kappa$ are
defined in the same way as in the analytic case.

 Now we compare both
definitions. For each $i=1,\ldots,r$, let $f_i:X\rightarrow
\A^{1}_{\C}$ be a separated reduced $\A^{1}_{\C}$-scheme of finite
type. Consider the strict localization $S$ of $\A_{\C}^{1}$ at
$0$, with closed point $s$, and generic point $\eta$. By abuse of
notation, we will also write $f_i:X\rightarrow S$ for the
restriction of $f_i$ to $X\times_{\A^{1}_{\C}}S$. Put $F=(f_i)$,
and let $R\psi^{an}_{F,I}(\mathcal{F})$ be the restriction of the
analytically defined complex of Alexander modules to $X_I=\{x\in
X_0\,|\,f_i(x)=0\  \mathrm{iff}\  i\in J\}$.

We can make use of an \'etale dictionary, allowing us to travel
from the algebraic to the analytic world, associating to a scheme
$X$ of finite type over $\C$ the analytic space $X(\C)$ of its
complex points, and to an \'etale sheaf $\mathcal{F}$ on $X_T$ a
sheaf $\mathcal{F}^{an}$ on $X_T(\C)$, with pleasant compatibility
properties - see, e.g. \cite{Frei}\cite{Milne}. Our Henselian ring
$S$ corresponds to the germ of a disc around the origin. The
universal cover of the one-dimensional torus, however, is not
algebraic in nature - this is exactly why we have to use limits in
the definition of the algebraic fundamental group \cite{Murre}. We
now explain how it is linked to a geometric generic point
$\bar{\eta}$ of $S$, see also \cite{SGA7b}.

Let $T_1$ be the torus $\A^{1}_{\C}\setminus \{0\}$, and
$\tilde{T}_1$ a universal cover. We can add a point $0$ to
$\tilde{T}_1$ to obtain a map $p:\tilde{\A}\rightarrow
\A^{1}_{\C}$ such that $p(0)=0$, $\tilde{T}_1$ is open in
$\tilde{\A}$, and the sets $p^{-1}(U)$, where $U$ runs over the
neighbourhoods of $0$ in $\A_{\C}^{1}$, form a fundamental system
of neighbourhoods of $0$ in $\tilde{\A}$. The cover $\tilde{T}_1$
determines an algebraic closure of the residue field $k(\eta)$ as
follows: $S$ is isomorphic to the ring of germs of holomorphic
functions on $\A^1_{\C}$ at $0$, that are algebraic over the field
of rational functions $K(\A_{\C}^{1})$. The field of germs at
$0\in \tilde{\A}$ of meromorphic functions on $\tilde{T}_1$ that
are algebraic over $K(\A^{1}_{\C})$, is an algebraic closure
$k(\bar{\eta})$ of $K(\A^{1}_{\C})$, hence of $k(\eta)$. The group
of deck transformations $\pi_1(T_1)$ acts on $k(\bar{\eta})$,
which yields a monomorphism
$$\Z=\pi_1(T_1)\rightarrow
\mathrm{Gal}(k(\bar{\eta})/k(\eta))=\hat{\Z}(1)(\C)\ (*).$$ Let
$\bar{\eta}$ be the geometric point corresponding to the inclusion
$k(\eta)\subset k(\bar{\eta})$.

 Our first comparison statement is the following:
\begin{prop}\label{comp1}
Let $\mathcal{G}$ be an object in $D_{c}^{b}(X_T,\Lambda)$, and
let $\mathcal{G}^{an}$ be the corresponding object in
$D_{c}^{b}(X_T(\C),\Lambda)$. Consider the constructible complex
$R\psi_{F,J}^{(t)}(\mathcal{G})^{an}$ on the stratum $X_J$ (or its
intersection with the inverse image of a small disc in $\C^{r}$).
It inherits a $\hat{\Z}(1)(\C)^{r}$ action, which induces, via the
morphism (*), a monodromy action of $\Z^{r}$. This latter object
is quasi-isomorphic to $R\psi_{F,J}^{an}(\mathcal{G}^{an})$.
\end{prop}
The proof is completely analogous to the proof in the univariate
case \cite{SGA7b}, and consists of a reduction to the normal
crossing-case for which we established formulas above.

Our next task is to compare the algebraic definition of
$R\psi_{F,I}^{(t)}$ with our definition of the Alexander complex.
 For each $i=1,\ldots,r$, let $f_i:X\rightarrow \A^{1}_{k}$ be
a morphism of $k$-varieties. Consider the strict
localization $S$ of $\A_{k}^{1}$ at $0$, with closed point $s$,
and generic point $\eta$. By abuse of notation, we will also write
$f_i:X\rightarrow S$ for the restriction of $f_i$ to
$X\times_{\A^{1}_{k}}S$. Put $F=(f_1,\ldots,f_r)$. Let $\xi$ be a
closed point of $X$ that is mapped to $s_J$; we might as well
assume that $J=\{1,\ldots,r\}$, and we drop the subscript $J$. Let
$\mu$ be the generic point of a strict localization of $X$ at
$\xi$, and let $\bar{\mu}$ be a geometric point localized at
$\mu$. We fix an isomorphism of fiber functors
$V_{F(\bar{\mu})}\cong V_{(1,\ldots,1)}$ on the category of tame
\'etale covers of $\G_{m,k}^{r}$. Consider an object $\mathcal{F}$
in $D_{c}^{b}(X_{\mathbb{G}^{r}_{m,k}},\bar{\Q}_l)$. We denote by
$R\psi^{(t),op}_{\mu,F}(\mathcal{F})$ the tame Alexander complex
at $\mu$, with opposite Gal$(\bar{\mu}/\mu)$-action.

Now we give $R\psi^{(t),op}_{\mu,F}(\mathcal{F})$ the structure of
a sheaf of modules on $\mathcal{C}(\G_{m,k}^{r})$, or rather: we
split up $R\psi^{(t),op}_{\mu,F}(\mathcal{F})$ in different direct
summands, in analogy with the Jordan decomposition of a linear
transformation, and each of these summands will get the structure
of a sheaf of modules on $\{\chi\}\times
\mathcal{C}(\G_{m,k}^{r})_l$, for some character $\chi$. For each
$\chi$ in $\mathcal{C}(\G_{m,k}^{r},\bar{\Q}_l^{\times})_f$ (see
\cite{GabLoe} for a definition), we denote by
$R\psi^{(t),op}_{\mu,F}(\mathcal{F})^{(\chi)}$ the biggest
subobject of $R\psi^{(t),op}_{\mu,F}(\mathcal{F})$ such that the
monodromy action, translated by $\chi^{-1}$, factors through
$\pi_{1}(\G_{m,k}^{r})_l$. Now it is clear how to give
$R\psi^{(t),op}_{\mu,F}(\mathcal{F})^{(\chi)}$ the structure of a
sheaf of modules on $\{\chi\}\times \mathcal{C}(\G_{m,k}^{r})$.

 For each closed point $\varphi$ of
$\{\chi\}\times \mathcal{C}(\mathbb{G}_{m,k}^{r})_l$, the stalk
$R\psi^{(t),op}_{\mu,F}(\mathcal{F})^{(\chi)}_{\varphi}$ is the
largest subobject of $R\psi^{(t),op}_{\mu,F}(\mathcal{F})$ on
which the action of each member $\gamma_j$ of a topological basis
$(\gamma_1,\ldots,\gamma_r)$ for $\pi_{1}(\GG_{m,k}^{r})^{t}$ is
of the form $\varphi(\gamma_j)+N_{\gamma_j}$, with $N_{\gamma_j}$
nilpotent.

\begin{prop}\label{monoco}
For each character $\chi$ in
$\mathcal{C}(\G_{m,k}^{r},\bar{\Q}_l^{\times})_f$, we get a
canonical isomorphism of sheaves of modules
$$\mathcal{A}^{q}_{F,\xi}|_{\{\chi\}\times
\mathcal{C}(\G_{m,k}^{r})_l}(\mathcal{F})\cong
R^{q-r}\psi^{(t),op}_{\mu,F}(\mathcal{F})^{(\chi)}\,.$$
\end{prop}
\begin{proof}
It suffices to prove the theorem when $\chi$ is trivial. We
replace $\mathcal{F}$ by a model in
$D_{c}^{b}(X_{\G_{m,k}^{r}},R)$, where $R$ is the ring of integers
of a finite extension of $\Q_l$. We denote this model again by
$\mathcal{F}$. Then the equality we want to prove becomes
$$\mathcal{A}^{q,R}_{F,\xi}(\mathcal{F})\cong
R^{q-r}\psi^{pro-l,op}_{\mu,F}(\mathcal{F})\otimes_{\Z_l}(\wedge_{\Z_l}^{r}\pi_1(\G_{m,k}^{r})_l)^{\vee}\,,$$
 where $^{\vee}$ denotes the dual $\Z_l$-module, and the
isomorphism is a canonical isomorphism of
$R[[\pi_1(\G_{m,k}^{r})_l]]$-modules.

We have a canonical isomorphism
$$\mathcal{A}^{q,R}_{F,\xi}(\mathcal{F})=R^{q}j_{*}(\mathcal{F}\otimes_{R}^{L}F^{*}L^{R}_{\G^{r}_{m,k}})_{\xi}
\cong
H^{q}(\mathrm{Gal}(\bar{\mu}/\mu),R\psi^{pro-l}_{\mu,F}(\mathcal{F}\otimes_R^{L}F^{*}L^{R}_{\G_{m,k}^{r}}))\,,$$
which may be rewritten as an isomorphism
$$\mathcal{A}^{q,R}_{F,\xi}(\mathcal{F})\cong
H^{q}(\pi_{1}(\G_{m,k}^{r})_l,R\psi^{pro-l}_{\mu,F}(\mathcal{F})\otimes_R^{L}L^{R}_{\G_{m,k}^{r}})\,,$$
where we view $L^{R}_{\G_{m,k}^{r}}$ as an $R$-module with
continuous $\pi_{1}(\G_{m,k}^{r})_l$-action. Now let us take a
look at the hypercohomological spectral sequence
\begin{eqnarray*}
E^{i,j}_{2}&=&H^{i}(\pi_1(\G_{m,k}^{r})_l,R^{j}\psi^{pro-l}_{\mu,F}(\mathcal{F})\otimes_R^{L}L^{R}_{\G_{m,k}^{r}})
\\&&\Rightarrow
H^{i+j}(\pi_1(\G_{m,k}^{r})_l,R\psi^{pro-l}_{\mu,F}(\mathcal{F})\otimes_R^{L}L^{R}_{\G_{m,k}^{r}})\,.\end{eqnarray*}
By \cite{GabLoe}, proposition 4.2.2.1, the spectral sequence
degenerates at $E_2$, and yields the desired isomorphism.
\end{proof}

%
%
What does this all mean for the support of the Alexander complex
$\mathcal{A}^{q}_{F,\xi}(\bar{\Q}_l)$ if $k=\C$? Consider the
universal cover
$$p:\C\rightarrow\C^{*}:z\mapsto e^{2\pi i z}\,,$$ and let
$\bar{\eta}$ be the corresponding geometric generic point of the
Henselization of $\A^{1}_{\C}$ at $0$. Fix an embedding of
$\bar{\Q}_l$ in $\C$. The choice of $\bar{\eta}$ yields a
distinguished topological basis $(\gamma_1,\ldots,\gamma_r)$ for
$\pi_1(\G_{m,\C}^{r})^t$, and the embedding of $\bar{\Q}_l$ in
$\C$ allows us to identify a closed point $\chi$ of
$\mathcal{C}(\G_{m,k}^{r})$ with a point $(e^{2\pi i\alpha_j})_j$
in Spec\,$\C[\Z]$, with $\alpha_j\in \Q$. This correspondence
identifies the support of $\mathcal{A}^{q}_{F,\xi}(\bar{\Q}_l)$
with the support of the stalk of the analytical complex
$\psi^{an}_{F}(\C)$ at $\xi$.

 To conclude, let us consider the behaviour of the
Alexander complex with respect to resolution of singularities. Let
$h:Y\rightarrow X$ be a simultaneous resolution for the $f_i$, and
let $\mathcal{F}$ be an object in $D_{c}^{b}(Y_{T},R)$. It follows
from the proper base change theorem that
$$Rh_{*}Rj_{*}(\mathcal {F}\otimes h^{*}F^{*}L^{R}_{T})\cong
Rj_{*}(Rh_{*}\mathcal{F}\otimes F^{*}L^{R}_{T})\,.$$ Suppose that
the resolution is tame, i.e. all multiplicities $N^{i}_{j}$ are
prime to the characteristic exponent of $k$. In this case, the
complex $Rj_{*}(\bar{\Q}_l\otimes h^{*}F^{*}L^{R}_{T})$ is lisse
on the strata $E_{I}^{o}$, by \cite{GabLoe}, lemma 4.3.2, and on
each stratum $E_{I}^{o}$, we have an isomorphism
$$R^{q}j_{*}(\bar{\Q}_l\otimes h^{*}F^{*}L^{R}_{T})\cong R^{0}j_{*}(\bar{\Q}_l\otimes h^{*}F^{*}L^{R}_{T})\otimes \wedge^{q}(M_{I}\otimes
\bar{\Q}_l(-1))\,$$ where $M_{I}$ is defined in the obvious way.
\section{A formula for the local zeta function}\label{form}
\noindent We recall, and introduce, some notation. Let $K$ be a
finite field extension of $\Q_p$, with ring of integers $R$, and
residue field $k\cong F_{q}$. We denote by $P$ the maximal ideal
of $R$. Let $\Phi$ be a residual, locally constant function on
$R^{n}$, and $F=(f_{1},\ldots,f_{r})$, with $f_i$ a polynomial in
$x=(x_1,\ldots,x_n)$ over $R$. Let
$\chi=(\chi_{1},\ldots,\chi_r)$, where $\chi_j$ is a
multiplicative character on $R^{\times}$ of order $d_{j}$, trivial
on $1+P$. Put $d=(d_1,\ldots,d_r)$. Denote by $N_i$ the vector
$(N_i^1,\ldots,N_i^r)$, and by $s$ the $r$-tuple of complex
variables $(s_1,\ldots,s_r)$. Let $\bar{.}$ denote reduction
modulo $P$. Let $(Y,h)=h:Y\rightarrow \A^{n}_K$ be a simultaneous
embedded resolution for $F$, and always assume tame good reduction
\cite{D}\cite{Denef7}\cite{Denef6}. Let $T$ be the index set of
the exceptional components, where $T_s$ indexes the components of strict
transform, and $T_e$ the exceptional divisors. The symbol $E_{I}$,
with $I\subset T$, means $\cap_{i\in I}E_i$, and $E_{I}^{o}$
stands for $E_{I}\setminus \cup_{j\in T\setminus I}E_j$.
Analogously, we write $\bar{E}_{I}$ instead of $\cap_{i\in
I}\bar{E}_i$, and $\bar{E}_{I}^{o}$ instead of
$\bar{E}_{I}\setminus \cup_{j\in T\setminus I}\bar{E}_j$.

We say that a multiplicity vector $N_i$ satisfies condition
$\gamma(N_i)$ if the character $\prod_{j=1}^{r}\chi_{j}^{N_i^{j}}$
is trivial. The condition $\gamma(N_i)$ replaces the condition
$d|N_i$ in the univariate case.
\begin{theorem}\label{formdenef}
Let $f_i\in R[x]$, $\bar{F}\neq 0$, and let $(Y,h)$ have good
reduction modulo $P$. Then
$$Z_{\Phi}(F,\chi;s)=q^{-n}\sum_{I\subset T, \gamma(N_i)\,if\, i\in
I}c_{I,\chi,\Phi}\prod_{i\in I}\frac{q-1}{q^{N_i s+\nu_i}-1}$$
where $c_{I,\chi,\Phi}=\sum_{a\in
\bar{E}_I^{o}(k)}\bar{\Phi}(a)\Omega_{\chi}(a)$ and
$\Omega_{\chi}(a)=\prod_{j=1}^{r}\chi_j(u_j)$, with $u_j$ as in
the univariate case \cite{Denef6} ($\chi_j(u_j)$ is ill-defined,
but their product isn't).
\end{theorem}

\begin{proof}
See the proof of Denef in \cite{Denef6}. His proof can be read as
if $s$, $\psi$ and the $N_i$ were vectors, replacing the condition
$d|N_i$ by $\gamma(N_i)$ to obtain the generalized formula.
\end{proof}
%
%

\noindent We write $c_{I,\chi}$ for $c_{I,\chi,\Phi}$ when $\Phi$
is the characteristic function of $R^{n}$, and $c_{I,\chi,0}$ when
$\Phi$ is the characteristic function of $P^{n}$.

These coefficients have a cohomological interpretation. Choose a
prime number $l$, not dividing $q$, such that all characters
$\chi_j$ take values in $\Q_l$, identifying $\Q_l$ with a subfield
of $\C$. Let $C$ be the torus $(\A_{F_{q}}\setminus 0)^{r}$,
define $\bar{Y}_C$ as $(\bar{F}\circ\bar{h})^{-1}(C)$ and let
$\alpha:\bar{Y}_C\rightarrow C$ be the induced map.

The \'etale cover
$$\mathrm{Spec}F_q[w,w^{-1},z]/(z_1^{d_1}-w_1,\ldots,z_r^{d_r}-w_r)$$ of
$C$, $z=(z_1,\ldots,z_r)$, $w=(w_1,\ldots,w_r)$, has Galois group
$\mu_d=\prod_{i=1}^{r} \mu_{d_{i}}$. Through the isomorphism
$\mu_d\cong\prod_{i}F_{q}^{\times}/(F_{q}^{\times})^{d_{i}}$,
$\psi=\prod_{i}\chi_i$ induces a character $\psi'$ on $\mu_{d}$,
which at its turn induces a
 continuous
representation of $\pi_1(C)$ in $\Q_l^{\times}$. We denote the
 locally constant $\Q_l$-sheaf associated to its inverse by $\mathcal{L}'_{\psi}$.
Let $\mathcal{L}_{\psi,\alpha}$ be
$\alpha^{*}\mathcal{L}'_{\psi}$, and $\mathcal{F}_{\psi}$ be
$\nu_{*}\mathcal{L}_{\psi,\alpha}$, where
$\nu:\bar{Y}_{C}\rightarrow \bar{Y}$ is the open immersion. The
sheaf is constructed in such a way that the action of the
geometric Frobenius corresponds to the character we want to study,
so we get a convenient Grothendieck trace formula \cite{SGA41/2}:
\begin{theorem}
The sheaf $\mathcal{F}_{\psi}$ is locally constant of rank 1 on
$U_{d}=\bar{Y}\setminus \cup_{\neg\gamma(N_i)}^{i\in T}\bar{E}_i$.
If $I\subset T$ and $\gamma(N_i)$ for all $i\in I$, then
$$c_{I,\chi} =\sum_{i}(-1)^{i}Tr(Fr,H_{c}^{i}(\bar{E}^{o}_{I}\otimes
_{F_{q}} F_{q}^{alg},\mathcal{F}_{\psi})),$$ and
$$c_{I,\chi,0} =\sum_{i}(-1)^{i}Tr(Fr,H_{c}^{i}(\bar{E}^{o}_{I}\cap \bar{h}^{-1}(0)\otimes
_{F_{q}} F_{q}^{alg},\mathcal{F}_{\psi}))\,$$ where $Fr$ is the
Frobenius endomorphism, and $F_{q}^{alg}$ denotes the algebraic
closure of $F_{q}$.
\end{theorem}
\begin{proof}
Cover $U_{d}$ with affine opens $V\subset U_d$, such that on each
$V$, we can write $\bar{f_j}\circ \bar{h}$ as $u_j \prod
y_i^{M^{j}_{i}}$, with $u_j$ nonvanishing on $V$, $j=1,\ldots, r$,
and $(y_j)$ a system of local parameters on $V$. Let $u$ be the
mapping from $V$ to the torus $C$ induced by the $u_j$. The
multiplicity vectors $M_i$ satisfy condition $\gamma$ for each
$i$, so it follows from the bimultiplicativity of the construction
of $\mathcal{F}$ that the sheaf
$\mathcal{G}=u^{*}\mathcal{L}_{\psi}$ on $V$ coincides with
$\mathcal{L}_{\psi,\alpha}$ on $V\cap \bar{Y}_{C}$. Since
$\mathcal{G}$ is, by purity, isomorphic to the direct image of its
restriction to $V\cap\bar{Y}_{C}$, we see that
$\mathcal{G}=\mathcal{F}_{\psi}|_{V}$.

Now let $I$ be a subset of $T$ so that $\gamma(N_i)$ for each
$i\in I$; thus $\bar{E}^{o}_{I}\subset U_d$. Let $a$ be a point in
$V$. The action of the geometric Frobenius $\varphi$ of $a$ on the
stalk of $\mathcal{F}_{\psi}$ at $a$, is induced by the monodromy
action on the stalk by the element of $\mu_{d}$ corresponding to
$\varphi$, and is given by $\Omega_{\chi}(a)$. Grothendieck's
trace formula yields the expressions for $c_{I,\chi}$ and
$c_{I,\chi,0}$.
\end{proof}
We now study the higher direct images of
$\mathcal{L}_{\psi,\alpha}$ under $\nu$.
\begin{lem}
(i) $R^{j}\nu_{*}\mathcal{L}_{\psi,\alpha}$ is zero outside
$U_{d}$ for all $j\geq 0$.

\noindent(ii) For $I\subset T$, let $\nu_{I}$ be the open
immersion $\nu_{I}:\bar{E}_{I}^{o}\rightarrow \bar{E}_{I}$. Then
$R^{j}\nu_{I*}(\mathcal{F}_{\psi}|_{\bar{E}_{I}^{o}})$ is zero
outside $\bar{E}_{I}\cap U_{d}$ for all $j\geq 0$.

\noindent(iii) Suppose $i\in T$ and $\gamma(N_i)$. Then
$(R^{j}\nu_{*}\mathcal{L}_{\psi,\alpha})|_{\bar{E}^{o}_{i}}$ is
zero for all $j\geq 2$, and canonically isomorhpic to
$\mathcal{F}_{\psi}(-1)|_{\bar{E}^{o}_{i}}$ for $j=1$.
\end{lem}
\begin{proof}
(i) Let $a$ be a closed point of $\bar{Y}\setminus U_d$. This
point belongs to $\bar{E}^{o}_{I}$ for some $I\subset T$, and $I$
contains an index $i_0$ such that $N_{i_{0}}$ does not satisfy
$\gamma$. We will prove that the stalk of
$R^{j}\nu_{*}\mathcal{L}_{\psi,\alpha}$ at $a$ is zero. This
assertion is local for the \'etale topology, which allows us to
replace $\bar{Y}$ by $A^{n}_{F_q}$, $\bar{Y}_{C}$ by $(\prod_{i\in
I}(A^{1}_{F_q}\setminus 0))\times A_{F_q}^{n-card\, I}$, $\alpha$
by $(\prod_{i\in I}y_{i}^{N^{j}_{i}})_j$, and $a$ by the origin.
Then we can write
$$\mathcal{L}_{\psi,\alpha}=\alpha^{*}\mathcal{L}_{\psi}=\boxtimes_{i\in
I}(\otimes_{j=1}^{r}\mathcal{L}_{\chi_{j}}^{\otimes N_i})\boxtimes
\mathcal{C}$$ where $\mathcal{C}$ denotes the constant sheaf
$\Q_l$ on $A_{F_q}^{n-card\, I}$, and $\mathcal{L}_{\chi_j}$
denotes the Kummer sheaf associated to $\chi_j$. Let $\mu$ be the
open immersion of $A^{1}_{F_q}\setminus 0$ in $A^{1}_{F_q}$. Since
$N_{i_0}$ does not satisfy $\gamma$, the stalk of
$\mu_{*}\otimes_{j=1}^{r}\mathcal{L}_{\chi_{j}}^{\otimes N_{i_0}}$
at the origin is zero (because we get a nontrivial monodromy
action). Local duality yields
$(R^{1}\mu_{*}\otimes_{j=1}^{r}\mathcal{L}_{\chi_{j}}^{\otimes
N_{i_0}})_0=0$. Hence
$(R^{j}\mu_{*}\otimes_{j=1}^{r}\mathcal{L}_{\chi_{j}}^{\otimes
N_{i_0}})_0=0$ for all $j\geq 0$. Now we can apply the K\"unneth
formula to the expression for $\mathcal{L}_{\psi,\alpha}$, to
obtain the desired result.

(ii) Because of (i), we may suppose that $\bar{E}^{o}_{I}\subset
U_d$. When $N_i$ satisfies $\gamma$, we get a trivial monodromy
action on $\otimes_{j=1}^{r}\mathcal{L}_{\chi_{j}}^{\otimes N_i}$,
so this sheaf is constant. Now one can use a local description
similar to the one in (i).

(iii) Apply purity results to the smooth pair consisting of
$\bar{E}_{i}^{o}$ and $\bar{Y}\setminus \cup_{k\neq i}\bar{E}_k$.
\end{proof}
The same arguments as in the univariate case \cite{Denef6} can be
used to prove the following result from the preceding lemma.
\begin{lem}
Let $E_{i_{0}}$ be proper, $\gamma(N_{i_{0}})$, and suppose that
$E_{i_{0}}$ intersects no $E_j$ with $\gamma(N_j)$, $j\neq i_0$.
Then
$$H^{i}_c(\bar{E}^{o}_{i_{0}}\otimes_{F_{q}}F_q^{alg},\mathcal{F}_{\psi})=0\
\mathrm{for}\ i\neq n-1\,.$$
\end{lem}
 So if in addition
$\chi(E_{i_{0}}^{o})=0$ (i.e. the topological Euler characteristic
of $E_{i_{0}}^{o}(\C)$ vanishes), we have no contribution in the
formula for $Z_{\Phi}(F,\chi;s)$, since $\chi(E_{i_{0}}^{o})$ and
$\chi_{c}(\bar{E}_{i}^{o}\otimes_{F_{q}}F_{q}^{alg},\Q_l)$ differ
by a nonzero factor, and this latter Euler characteristic with
compact support at its turn equals
$\chi_{c}(\bar{E}_{i}^{o}\otimes_{F_{q}}F_{q}^{alg},\mathcal{F}_{\psi})$,
because $\mathcal{F}_{\psi}|_{\bar{E}_{i}^{o}}$ has tame
ramification.
 This
allows us to prove:
\begin{theorem}\label{unmask}
Suppose $n=2$. If $H: \sum_k c_k s_k+1=0$ is a polar hyperplane of
$Z$, then there exists an index $j$ in $T$, with
$|E_j/E^{o}_j|\geq 3$ or $j\in T_s$, and $c_k=N_j^{k}/\nu_j$ for
each $k$.
\end{theorem}
\begin{proof}
This proof generalizes the one in \cite{DenefBour}. Some notation:
fix $j$ in $T_e$, and let $a_t$, $t\in J$, be the geometric points
of $E_j\setminus E_j^{o}$. For $t\in J$, let $N^{k}_t,\nu_t$
denote the numerical data associated to the unique divisor
intersecting $E_j$ transversally in $a_t$, and put
$\alpha_t^{(k)}=\nu_t-N^{k}_t\nu_j/N^k_j$ if $N_j^k$ differs from
zero. It is well-known \cite{Loepadic}\cite{Veysnum} that
$$\sum_{t\in J}(\alpha_t^{(k)}-1)=-2,\ \mathrm{and}\ \sum_{t\in
J}N_t^{k}\equiv 0\ mod\,N_j^{k}\,\quad(*)$$ for each $k$ for which
the identities are defined.

 Let $S$ be the set of all divisors $E_j$ inducing the polar
hyperplane $H$. Suppose that each of its members satisfies
$|E_j\setminus E^{o}_j|\leq 2$ and $j\notin T_s$. It follows from
the identities above that different elements of $S$ cannot meet,
since otherwise we could build an infinite tree of divisors in
$S$, which is absurd. Suppose that $E_j$ belongs to $S$ and
satisfies $\gamma(N_j)$. Suppose that $|E_j\setminus
E^{o}_j|=\{a_1,a_2\}$; the other case is similar. Either both
$N_1$ and $N_2$ satisfy $\gamma$, or none of them does; in the
latter case $E_j$ does not contribute to $Z$ because of the
previous lemma, in the former case the sheaf $\mathcal{F}_\psi$
will be geometrically constant on $\bar{E}_j$, and it follows from
the trace formula and the formula for $Z$, combined with the
information about the $\alpha_t^{(k)}$, that $E_j$ does not
contribute to the polar locus of $Z$.
\end{proof}

\section{Degree of local zeta functions and
monodromy}\label{degree} \noindent In this section, we establish a
formula for the limit of the zeta function for $s\to - \infty$,
and we prove the following theorem, which is a generalization of a
theorem in \cite{Denef5}:
\begin{theorem}\label{theoremdegree}
Let $F=(f_1,\ldots,f_r)$ be defined over a number field $K\subset
\C$. For almost all completions of $K$, the degree of the local
zeta function $Z_{0}(f,\chi;s)$ will be strictly less then zero if
$\chi\notin \mathrm{Supp}\mathcal{A}^{q}_{F,0}(\bar{\Q}_l)$ for
all $q$.
\end{theorem}
Here the statement that the degree of $Z_{0}(f,\chi;s)$ is
strictly negative, means that the well-defined limit $\lim_{s\to
-\infty}Z_{0}(f,\chi;s)$ vanishes.
\begin{proof}
\noindent Let $\bar{\eta}$ be a geometric generic point of
$\A^{r}_{F_{q}}$, defining an algebraic closure $F_{q}^{alg}$ of
$F_{q}$. Let $S_0$ be the Henselization at $0$ of
 $\A^{r}_{F_{q}}$, and denote by
 $\eta_0$ its generic point. Put
$G_0=\mathrm{Gal}(\bar{\eta}/\eta_0)$.


Let $\sigma$ be an element of $G_0$ which induces the geometric
Frobenius $\varphi$ on $F_{q}^{alg}$. Then, using the notation
from the previous section,

\begin{eqnarray*}
&&\sum_a(-1)^{a}Tr(\sigma,\mathcal{A}^{a}_{\bar{F},0}(\bar{\Q}_l)_{\chi})
\\&=&\sum_{a,b}(-1)^{a+b}Tr(\sigma,H^{a}(\bar{h}^{-1}(0)\otimes
F_{q}^{alg},R^{b}\nu_{*}((\bar{F}\circ
\bar{h})^{*}L^{\Z_l}_{T}\otimes (\bar{F}\circ
\bar{h})^{*}\mathcal{L}_{\chi}))
\end{eqnarray*}
(By proper base change and Leray's direct image spectral sequence,
we can compute the trace on the exceptional locus of the embedded
resolution $\bar{h}$.)
\begin{eqnarray*}
&=&\sum_{I}\sum_{a,b}(-1)^{a+b}Tr(\sigma,H^{a}_{c}((\bar{E}^{o}_{I}\cap\bar{h}^{-1}(0))\otimes
F_{q}^{alg}, \\[-10pt]&&\qquad \qquad \qquad R^{0}\nu_{*}(
(\bar{F}\circ \bar{h})^{*}L^{\Z_l}_{T}\otimes (\bar{F}\circ
\bar{h})^{*}\mathcal{L}_{\chi})\otimes \wedge^{b}(M_{I}\otimes
\bar{\Q}_l(-1)))
\end{eqnarray*}
(By the explicit computations in the normal crossing case in
Section \ref{trans})
\begin{eqnarray*}
&=&\sum_{I}\sum_{a,b}(-1)^{a+b}Tr(\varphi,H^{a}_{c}((\bar{E}^{o}_{I}\cap\bar{h}^{-1}(0))\otimes
F_{q}^{alg},\mathcal{F}_{\chi}))Tr(\varphi,\wedge^{b}(M_{I}\otimes
\bar{\Q}_l(-1)))
\\&=&q^{-n}(1-q)^{\lambda}\lim_{s\to -\infty}Z_{0}(F,\chi,s)
\end{eqnarray*}
where $\lambda$ is the number of linearly independent multiplicity
vectors $(N^{j}_{i})_{i\in I}$, $j=1,\ldots,r$.
\end{proof}

In the motivic setting, the limit of the multivariate zeta
function yields a motivic incarnation of the Alexander complex
\cite{Guibert}. In the case $r=1$, the motivic Milnor fiber has
the same mixed Hodge structure as the geometric Milnor fiber, and
 the Hodge spectrum of $F=f_1$ can be recovered from this motivic object
\cite{DL5}\cite{DL3}.
\section{Proof of the monodromy and holomorphy conjectures for
curves}\label{proofs}

Now we are ready to prove Conjectures \ref{gmonoconjecture} and
\ref{gholoconjecture} when $n=2$ (see \cite{Rod} and \cite{Veys1}
for the univariate case $r=1$).

\begin{proof}[Proof of the Generalized Monodromy Conjecture \ref{gmonoconjecture} for $n=2$]
Polar hyperplanes induced by components of the strict transform
satisfy the conjecture, as can be seen by looking at the Alexander
polynomial of the stalk of the Alexander complex in a smooth point
of the zero locus. Let $E_j$ be an exceptional divisor
intersecting at least three other divisors and lying above a point
$x$ of $X$. Let $N'$ be an integer $r$-tuple, and denote by
$\mathcal{I}$ the set of exceptional components lying over $x$,
with a multiplicity vector which is an integer multiple of $N'$.
We say that $E_i$ satisfies $\delta$ if it belongs to
$\mathcal{I}$. Because of Sabbah's A'Campo formula, and the
previous arguments, it suffices to prove that the expression
$w_{\mathcal{I}}:=\sum_{E_{i}\in \mathcal{I}} \chi(E_i^{o})$ is
strictly less than zero, for each $N'$ such that $E_j$ satisfies
$\delta$. We may suppose that none of the components of the strict
transform satisfies $\delta$, since in the other case, it again
suffices to look at the Alexander complex in a smooth point. To
prove the monodromy conjecture, we use the inequality
$w_{\mathcal{I}}<0$ in the case where $N'$ equals $N_j$, divided
by the greatest common divisor of $\nu_j$ and the entries of
$N_j$.

Let $\mathcal{N}$ be the union over all $E_i$ satisfying $\delta$.
It follows from the relations (*) in the proof of Theorem
\ref{unmask}, that each connected component of $\mathcal{N}$
contains a prime divisor, intersecting at least two divisors $E_k$
not satisfying $\delta$ (we suppose $f_j(x)=0$ for some $j$). Let
$\mathcal{M}$ be the set of divisors in the connected component of
$\mathcal{N}$ containing $E_j$; it suffices to prove that
$w_{\mathcal{M}}<0$. It is known that $w_\mathcal{M}=0$ implies
that $\mathcal{M}$ is a tree intersecting the union of the other
exceptional components exactly twice \cite{Rod}. These two
intersection points have to lie on the same divisor $E_1$ in
$\mathcal{M}$. Suppose that $E_1=E_j$. There is still another
divisor $E_2$ in $\mathcal{M}$ intersecting $E_j$. We claim that
we can find an index $k$ in $\{1,\ldots,r\}$ such that
$\alpha^{(k)}_t<1$ for all $E_t$ intersecting $E_1$. Since
$\alpha^{(k)}_2\in \Z$, (*) yields $\alpha^{(k)}_2=0$. But $E_2$
has to intersect another divisor in $\mathcal{M}$, and repeatedly
applying (*) produces an infinite tree. If $E_1\neq E_j$, there
are at least three divisors in $\mathcal{M}$ intersecting $E_j$;
since only one component intersects a divisor not belonging to
$\mathcal{M}$, the same arguments can be used.

We still have to prove our claim. We claim that we can always find
an appropriate index $k$, unless we have done some redundant
blow-ups during the resolution process. It is known that this
property holds in the univariate case $r=1$. But all the divisors
meeting $h^{-1}(x)$ appear as well in a minimal resolution for
$\prod_{f_j(x)=0}f_j$, so the claim follows from
$\sum_{i=1}^{m}x_i/\sum_{i=1}^{m}y_i\leq \max_i\{x_i/y_i\}$ for
$x_i, y_i\in \N_0$.
\end{proof}

\begin{proof}[Proof of the Generalized Holomorphy Conjecture
\ref{gholoconjecture} for $n=2$] If $Z$ is not holomorphic, there
must exist an index $j$ in $T$, such that the corresponding
divisor $E_j$ induces a pole in the formula for the zeta function.
This means that we can choose $j$ such that $\gamma(N_j)$ is
satisfied, and $j\in T_s$ or $|E_j\setminus E_j^{o}|\geq 3$. If
$j\in T_s$, it is clear from Sabbah's A'Campo formula that the
holomorphy conjecture is satisfied; so let us assume that
$|E_j\setminus E_j^{o}|\geq 3$, and that $E_j$ is an exceptional
divisor above $x$. We've proven above that that the expression
$w_I:=\sum_{E_{i}\in I} \chi(E_i^{o})$ is strictly less than zero
when we take the sum over all exceptional divisors lying over $x$
for which the multiplicity vector is a multiple of some integer
vector $N'$, whenever $E_j$ belongs to $I$.

Taking $N'$ equal to $N_j$, defining $g$ to be the greatest common
divisor of the entries of $N_j$, and putting $N'_j=N_j/g$, it
follows from Sabbah's A'Campo formula that the locus of
$T^{N'_j}=\xi_{g}$ is contained in the support of the Alexander
complex at $x$, where $\xi_g$ is a primitive $g$-th root. Hence,
it suffices to find a character $\psi$ so that $\psi^{N'_j}$ is a
primitive $g$-th root, and $\chi=\psi^{k}\chi'$, with
$\chi'^{N'_j}=1$ and $k\in \N$. Since $\gamma(N_j)$ is satisfied,
we can take $\chi'$ to be $\chi^{g}$, and furthermore, if $\chi$
is non-trivial, we can find $k\in \N$ so that
$\chi^{N_j(1-g)/k}=1$, but $\chi^{g'N'_j(1-g)/k}\neq 1$ unless $g$
divides $g'$. Now put $\psi=\chi^{(1-g)/k}$.
\end{proof}
\section*{Acknowledgements}
The author would like to thank Fran\c{c}ois Loeser for his
clarifications and comments, and for the inspiring discussions.
 \nocite{Loe2} \nocite{prox1}
\nocite{prox2}

\bibliographystyle{amsplain}
\bibliography{wanbib,wanbib2}
\end{document}